\documentclass[12pt]{article}
\usepackage{graphicx}
\usepackage{amssymb}
\newcounter{fig}

\newtheorem{theorem}{Theorem}

\title{Minimizing the footprint of your laptop \\(on your bedside table)}
\author{Burkard Polster}
\begin{document}

\maketitle

\section{Introduction} I often work on my laptop in bed.  When needed, I park the laptop on the bedside table, where the computer has to share the small available space with a lamp, books, notes, and heaven knows what else. It often gets quite squeezy.

Being regularly faced with this tricky situation, it finally occurred to me to determine once and for all how to place the laptop on the bedside table so that its ``footprint'' - the area in which it touches the bedside table - is minimal. In this note I give the solution of this problem, using some  very pretty elementary mathematics. 

\section{Mathematical laptops and bedside tables}

We assume that both the laptop and the bedside table are rectangular, and we will refer to these rectangles as the {\em laptop} and the {\em table}. We further assume that the center of gravity of the laptop is its midpoint. Finally, without loss of generality we may  assume that the laptop is 1 unit wide.\footnote{That is, the shorter side of the laptop is 1 unit in length. And, if the laptop is square then it is a unit square. Yes, yes, only a mathematician would consider  the possibility of a square laptop, but bear with me. As will become clear, considering square laptops provides an elegant key to our problem.} 

We are considering all placements of the laptop such that it  will not topple off the table; these are exactly the placements for which the midpoint of the laptop is also a point of the table. We are then interested in determining for which of these  placements the {\em footprint} of the laptop is of minimal area; here, the footprint is  the common region of the laptop and the table. 

\begin{figure}[h]
\centerline{\includegraphics{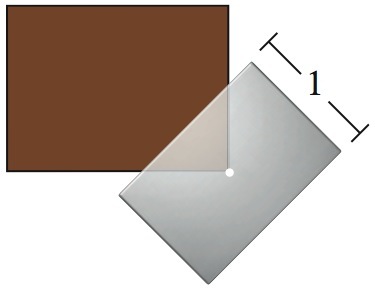}}
\caption[f1]{No (stable) placement of your laptop on a bedside table has a smaller footprint.\label{laptop1}}
\end{figure}

In all reasonable circumstances, the optimal  answer to this problem will always resemble the arrangement in Figure~\ref{laptop1}.\footnote{``Reasonable circumstances'' means in reference to laptops and tables of relative dimensions close to those of the real items. In the nitty gritty of this note we'll specify  the exact scope of our solution, and also what happens in some unrealistic but nevertheless mathematically interesting scenarios.}  {\em This optimal placement is characterized by the fact that the midpoint of the laptop coincides with one of the corners of the table and the footprint is an isoceles right triangle.}

The proof is divided into two parts. First, we consider those placements for which the midpoint of the laptop coincides with one of the corners of the table: we prove that among such placements our special placement has smallest footprint area. Then, we extend  our argument, proving that any placement for which the laptop midpoint is not a table corner must have a greater  footprint area.

\section{Balancing on a corner}

\begin{figure}[h]
\centerline{\includegraphics{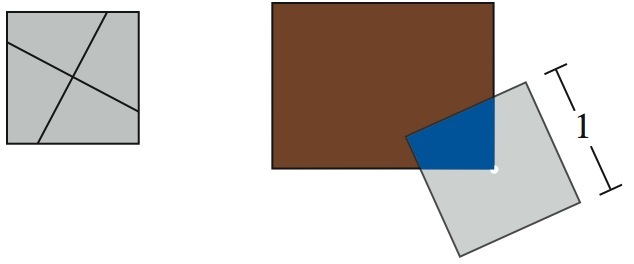}}
\caption[f1]{A square laptop with midpoint at a  corner will have a footprint area of 1/4. \label{square}}
\end{figure}

We begin by considering a right-angled cross through the center of a square, as illustrated in the left diagram in Figure~\ref{square}. Whatever its  orientation, the cross cuts the square into four congruent pieces. This shows that if we place a unit square laptop on the corner of a sufficiently large  table, its footprint will always have area 1/4, no matter how the square is oriented; see the diagram on the right. 

\pagebreak

\begin{figure}[h]
\centerline{\includegraphics{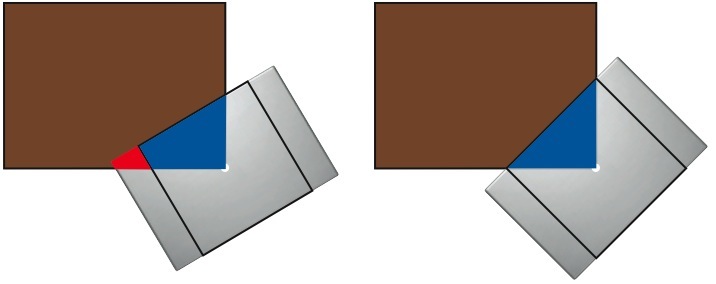}}
\caption[f1]{The blue regions have the same area, and so the right footprint is smaller. \label{overlap}}
\end{figure}

Next, consider a non-square laptop with its midpoint on the corner of a large table, as in Figure~\ref{overlap}. We regard the footprint as consisting of a blue part and a red part, as shown. Our previous argument shows that, as we rotate the laptop,  the  area of the blue part stays constant. On the other hand,  the red part only vanishes in the special position shown on the right. We conclude that this symmetric placement of the laptop uniquely provides the footprint of least area.

\begin{figure}[h]
\centerline{\includegraphics{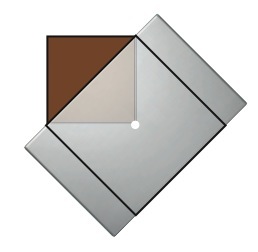}}
\caption[f1]{Our corner argument works for relatively small tables. \label{smallsquare}}
\end{figure}

These arguments required  that the  table be sufficiently large. How large? The arguments work as long as the rotated square never
pokes off another side of the table. So, since the short side of the laptop has length 1, we only require that
the shortest side of the  table be at least of length $1/\sqrt 2$; see Figure~\ref{smallsquare}. 
\pagebreak

\section{In the corner is best}

We now want to convince ourselves that the minimal footprint must occur for one of these special placements over a table corner. 

\begin{figure}[h]
\centerline{\includegraphics{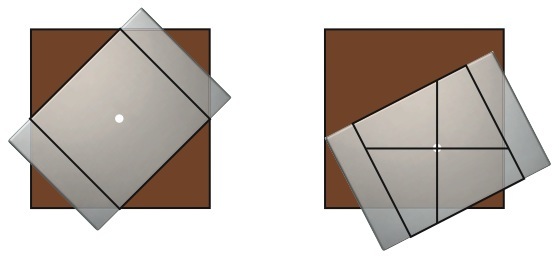}}
\caption[f1]{If the table contains the square highlighted on the left, then at least one of the quarters of the square on the right is contained in the footprint of the laptop. \label{corner_best}}
\end{figure}

We start with a table that is at least as wide as the diagonal of the square inscribed in our laptop; see the left diagram in Figure~\ref{corner_best}. Place the laptop anywhere on the table.  Now  consider a cross in the middle of the laptop square, and 
with arms parallel to the table sides. 

 As we saw above, the cross cuts the square into four congruent pieces. Furthermore,   wherever the laptop is placed 
and however it is oriented, at least one of these congruent pieces will be part of the footprint:  this is a consequence
of  our assumption on the table size. Finally,
 unless the midpoint is over a corner of the table, this quarter-square region clearly cannot be the full footprint.

Putting everything  together, we can therefore guarantee that our symmetric corner arrangement  is optimal if the table is at least as large as the square table in Figure~\ref{corner_best}. This square table has side length $\sqrt 2$. 

By refining the arguments above, we now want to show that our solution holds for any table that is at least 1 unit wide.  Since our laptop is also 1 unit wide, this probably takes care of most real life laptop balancing problems. 

\begin{figure}[h]
\centerline{\includegraphics{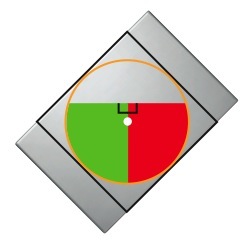}}
\caption[f1]{Both the red and the green regions have the critical area of  1/4.  \label{red_square}}
\end{figure}

Begin with a circle inscribed in the laptop square, and with the red and the green regions within, as in Figure~\ref{red_square}. The
regions are mirror images, and are arranged to each have area 1/4. Note that if the laptop is rotated around its midpoint,
either fixed region remains within the laptop. 

\begin{figure}[h]
\centerline{\includegraphics{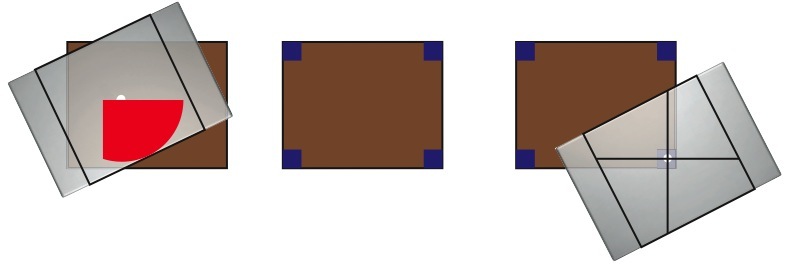}}
\caption[f1]{The footprint area is at least 1/4, for the laptop midpoint in either the blue or  brown region. \label{thin1}}
\end{figure}

Now place the laptop on the table with some orientation.  Suppose that the laptop  footprint contains a red or a  green region,
or such a region rotated by 90, 180, or 270 degrees; see Figure~\ref{thin1}. 
Then it is immediate that the footprint area for the laptop in that position is greater than~1/4. 

In fact the footprint may not contain such a region. However, this will be the case unless the laptop midpoint is close to a table corner,
in one of the little blue squares pictured in Figure~\ref{thin1}.
 On the other hand, if the midpoint is in a blue  square
  then the footprint will contain one of the original quarter-squares of area 1/4; see the diagram on the right side of Figure~\ref{thin1}.  

At this point we summarize what we have discovered so far.

\begin{theorem} Consider a laptop that is 1 unit wide and a table that is at least 1 unit wide. If the laptop is not a square, then the placement of the laptop on the table that gives the smallest footprint is shown in Figure~\ref{laptop1}. If the laptop is a square, then  the  minimal area footprints are for placements for which the midpoint of the laptop coincides with a corner of the table. 
\end{theorem}

\section{Odds and ends}

\begin{figure}[h]
\centerline{\includegraphics{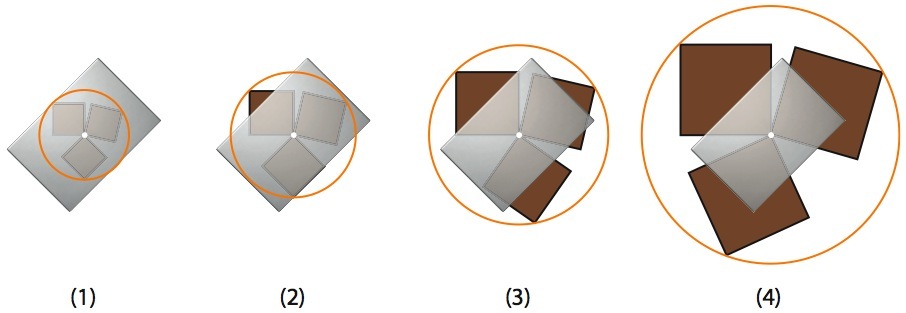}}
\caption[f1]{Balancing a laptop on the corners of squares of different sizes.  \label{thin3}}
\end{figure}

What  if you are the unlucky owner of a really small bedside table? First of all,  it is usually not 
difficult to determine the best placement for a specific laptop/table combination. To get a feel for this, and for
what to expect in general, consider Figure~\ref{thin3}, where we balance a laptop on square tables of different sizes. 

Here are some simple observations, applicable both to square and non-square tables:

\begin{enumerate}
\item If your table is really tiny, the footprint will always be the whole table, no matter where the laptop is placed.
This will be the case if  the table  diagonal is no longer than half the width of the laptop.
\item Suppose that the (square or non-square) table diagonal is just a little bit longer than half the width of the laptop, ensuring that if part of the table sticks out from underneath the laptop, then this part is a triangle cut off one of the corners of the table. 

In this case, if a table corner is sticking out and if the laptop midpoint is not at the opposite corner, then it is easy to see that simply translating the laptop to this opposite corner will lower the footprint area. Consequently, the minimal footprint will correspond to one of these special placements. For a square laptop, the minimal footprint occurs when the protruding triangle is isosceles. This is not terribly surprising. What may be surprising is that this is not at all obvious to prove; a descent into the land of nitty gritty seems unavoidable. 

For non-square tables, the optimal placement will not necessarily correspond to an isosceles table corner  sticking out. To see this, consider a very thin table. Then the only way a corner  can stick out  is if the table diagonal is almost perpendicular to the long side of the laptop. This precludes an isosceles triangle part of the table sticking out; see Figure~\ref{thin4}.

\begin{figure}[h]
\centerline{\includegraphics{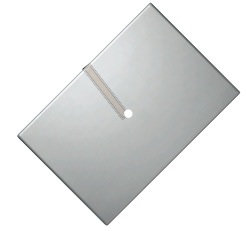}}
\caption[f1]{Balancing the laptop on a very thin table with a tiny piece of the table showing.  \label{thin4}}
\end{figure}
\item From here on things get even more complicated:  all the problems that we mentioned for the last  scenario plus many-sided, odd-shaped footprints, no easy way to see why the best placement should be among the placements for which the midpoint of the laptop is one of the corners, etc.
\item From here on our theorem applies.
\end{enumerate}

We end this note with some challenges for the interested reader (in likely order of difficulty):
\begin{itemize}
\item Extend our theorem to include all tables that are at least $1/\sqrt 2$ wide. (The table shown in Figure~\ref{smallsquare} has these dimensions.)
\item Turn Scenario 2 discussed above into a theorem (are there pretty proofs?).
\item Prove the Ultimate Laptop Balancing Theorem, that includes everything that your lazy author did and did not cover in this note: arbitrary location of the center of gravity, starshaped laptops and jellyfish-shaped tables, higher-dimensional tables and laptops, etc.
\end{itemize}

Have Fun, and Good Luck! 

\pagebreak

\noindent Burkard Polster\\
School of Mathematical Sciences\\
Monash University, Victoria 3800\\
Australia\\
e-mail: Burkard.Polster@sci.monash.edu.au\\
web: www.qedcat.com

\end{document}